\documentclass[11pt]{article}
\usepackage{amsfonts}
\usepackage{latexsym, amssymb, amsmath, amscd, amsfonts, epsfig, graphicx, colordvi}
\usepackage{ifpdf}
\usepackage{graphicx}
\newtheorem{thm}{Theorem}[section]

\newtheorem{defi}[thm]{Definition}

\def\pf{\noindent{\it Proof.} }
\def\qed{\nopagebreak\hfill{\rule{4pt}{7pt}}
\medbreak}

\def\qed{\nopagebreak\hfill{\rule{4pt}{7pt}}
\medbreak}

\makeatletter
\def\ExtendSymbol#1#2#3#4#5{\ext@arrow 0099{\arrowfill@#1#2#3}{#4}{#5}}
\makeatother

\title{Anti-lecture Hall Compositions and Overpartitions}
\author{William Y.C. Chen\raisebox{5pt}{\scriptsize 1},
Doris D.M. Sang\raisebox{5pt}{\scriptsize 2}, and Diane Y.H. Shi\raisebox{5pt}{\scriptsize 3}}
\date{Center for Combinatorics, LPMC-TJKLC\\
 Nankai University\\
Tianjin 300071, P.R. China \\
\vspace{15pt}
\raisebox{5pt}{\scriptsize 1\,}chen@nankai.edu.cn,
\raisebox{5pt}{\scriptsize 2\,}sdm@cfc.nankai.edu.cn, \raisebox{5pt}{\scriptsize 3\,}yahuishi@gmail.com}

\begin{document}
%--------------------------------------------------------------------
\maketitle
%--------------------------------------------------------------------
\noindent {\bf Abstract.}
 We show that the number of  anti-lecture
hall compositions of $n$ with the first entry not exceeding $k-2$ equals the number of overpartitions of $n$
with non-overlined parts not congruent to $0,\pm 1$ modulo $k$. This identity  can be considered as a refined version of the
anti-lecture hall theorem of Corteel and Savage.
To prove this result, we find two Rogers-Ramanujan type identities for overpartition which are analogous to the Rogers-Ramanjan type identities due to Andrews.
When $k$ is odd, we give an alternative proof by using a generalized Rogers-Ramanujan identity due to Andrews, a bijection of Corteel and Savage and a refined version of a bijection also due to Corteel and Savage.


%--------------------------------------------------------------------

\section{ Introduction}

The objective of this paper is to establish a
connection between anti-lecture hall compositions
with an upper bound on the first entry and overpartitions under a congruence condition
on non-overlined parts.

In \cite{cor02}, Corteel and Savage introduced
the notion of anti-lecture hall compositions and obtained a formula for the generating function by constructing a bijection.
An anti-lecture hall composition of length $k$ is defined to be an integer sequence $\lambda=(\lambda_1,\lambda_2,\ldots,\lambda_k)$ such that
\[\frac{\lambda_1}{1}\geq\frac{\lambda_2}{2}\geq\cdots\geq\frac{\lambda_{k-1}}{k-1}\geq\frac{\lambda_k}{k}\geq0.\]
The set of anti-lecture hall compositions of length $k$ is denoted by $A_k$.
Corteel and Savage have shown that
\begin{equation}\label{anti}
\sum_{\lambda\in A_k}q^{|\lambda|}=\prod_{i=1}^k\frac{1+q^i}{1-q^{i+1}}.
\end{equation}
Let $A$ denote
the set of anti-lecture hall compositions. Since any anti-lecture hall composition can be written as
 an infinite vector ending with zeros, we have $A= A_\infty$ and
\begin{equation}\label{anti2}
\sum_{\lambda\in A}q^{|\lambda|}=\prod_{i=1}^{\infty}\frac{1+q^i}{1-q^{i+1}}.
\end{equation}

In view of the above generating function, one sees that
anti-lecture hall compositions are related to overpartitions. An overpartition of $n$ is defined by  a non-increasing sequence of natural numbers whose sum is $n$
in which the first occurrence of a number may be overlined, see,  Corteel and Lovejoy \cite{cor04}. In the language of overpartitions, the right side of (\ref{anti2})
is the generating function for  overpartitions of $n$ with the non-overlined parts larger than $1$.

The main result of this paper is the following
refinement of the anti-lecture hall theorem of
Corteel and Savage \cite{cor02}:

\begin{thm}\label{m-thm-1} For $k\geq 3$,
\begin{equation}
 \sum_{\lambda_1\leq k-2, \lambda\in A}q^{|\lambda|}=\frac{(-q;q)_{\infty}}{(q;q)_{\infty}}
(q;q^{k})_{\infty}(q^{k-1};q^{k})_{\infty}(q^{k};q^{k})_{\infty}.
\end{equation}
\end{thm}

We shall make a connection
between  anti-lecture hall
compositions and the  overpartions with congruence restrictions.
Let $F_k(n)$ be the set of anti-lecture hall compositions  $\lambda=(\lambda_1,\lambda_2,\ldots)$ of $n$ such that
$\lambda_1 \leq k$.
Let $H_k(n)$ be the set of overparitions of $n$ for which the non-overlined parts are not
congruent to $0,\pm 1$ modulo $k$.
 Therefore, Theorem \ref{m-thm-1} can be restated
 as the following equivalent form.

  \begin{thm}\label{main}
For $k\geq 3$ and any positive integer $n$, we have
\begin{equation}
 |F_{k-2}(n)|=|H_k(n)|.
\end{equation}
\end{thm}

To prove the main result, we need to compute the
generating functions of the anti-lecture
hall compositions $\lambda$ with $\lambda_1\leq k$,
depending on the parity of $k$.
Then we shall show that these two generating functions of the anti-lecture hall compositions
 in $F_{2k-2}(n)$ and $F_{2k-3}(n)$
are equal to the generating functions of overpartitions in $H_{2k}(n)$ and $H_{2k-1}(n)$ respectively.
To this end, we  establish two Rogers-Ramanujan type identities (\ref{id1}) and (\ref{id2}) for overpartitions which are analogous to
  the following Rogers-Ramanujan type identity obtained by Andrews \cite{and74, and75}:
  \begin{equation}\label{generr}
\sum_{N_1\geq N_2 \geq \cdots \geq N_{k-1}\geq 0}
\frac{q^{N_1^2+N_2^2+\cdots+N_{k-1}^2+N_a+\cdots+N_{k-1}}}{(q)_{n_1}(q)_{n_2}\cdots(q)_{n_{k-1}}}
=\frac{(q^{a}; q^{2k+1})_{\infty}(q^{2k+1-a}; q^{2k+1})_{\infty}(q^{2k+1}; q^{2k+1})_{\infty}}{(q;q)_{\infty}}
\end{equation}
where $n_i=N_i-N_{i+1}$ and $1\leq a \leq k$. For $k=2$ and $a=1,2$, (\ref{generr}) implies the classical
Rogers-Ramanujan identities \cite{gor61}:
\begin{equation}\label{rr1}
\sum_{n=0}^{\infty}\frac{q^{n^2}}{(q)_n}=\prod_{n=0}^{\infty}(1-q^{5n+1})^{-1}(1-q^{5n+4})^{-1}
\end{equation}

\begin{equation}\label{rr2}
\sum_{n=0}^{\infty}\frac{q^{n^2+n}}{(q)_n}=\prod_{n=0}^{\infty}(1-q^{5n+2})^{-1}(1-q^{5n+3})^{-1}.
\end{equation}
It is worth mentioning that  Andrews' multiple series transformation \cite{and75} can be employed to derive the overpartition analogues of (\ref{generr}).

 When the upper bound $k$ is even, the
  weighted counting anti-lecture hall compositions leads to the left hand side
  of the first Rogers-Ramanujan type identity (\ref{id1}), whereas the generating
  function for the number of overpartitions equals the right hand side of
 the first Rogers-Ramanujan type identity (\ref{id1}). The case when $k$ is odd can be dealt with in  the same way.

When $k$ is odd, we  provide an alternative proof based on a refined version of  a bijection of Corteel and Savage
\cite{cor02}, a bijection of Corteel and Savage in the original form for the anti-lecture hall
theorem,  and a generalized Rogers-Ramanujan identity (\ref{generr}) of Andrews.

This paper is organized as follows: In Section 2, we give two Rogers-Ramanujan type identities for overpartitions.
Section 3 is concerned with the case of an even upper bound $k$. Two  proofs
for the case of an odd upper bound will be presented in Section 4.

\section{Rogers-Ramanujan type identities for overpartitions}

In this section, we  give two Rogers-Ramanujan type identities (\ref{id1}) and (\ref{id2})
for overpartitions. It can be seen that the right side of (\ref{id1}) is the generating function
for overpartitions in $H_{2k}(n)$. In the next section we shall show that the left side of (\ref{id1}) equals
the generating function for anti-lecture hall compositions in $F_{2k-2}(n)$.
Similarly, the right side of (\ref{id2}) equals the generating function for overpartitions in $H_{2k-1}(n)$. In
Section 4 we shall show that the left side of (\ref{id2}) equals the generating function for anti-lecture hall
 compositions in $F_{2k-3}(n)$.

Let us recall Andrews' multiple series transformation \cite{and75}:
\begin{align}\label{andtr}
&_{2k+4}\phi_{2k+3}
\left[{ {a,q\sqrt{a},-q\sqrt{a},b_1,c_1,b_2,c_2,\ldots,b_k,c_k,q^{-N};q,\frac{a^k q^{k+N}}{b_1\cdots b_kc_1\cdots c_k}}
 \atop  {\sqrt{a},-\sqrt{a},aq/b_1,aq/c_1,aq/b_2,aq/c_2,\ldots,aq/b_k,aq/c_k,aq^{N+1}}}\right]\nonumber\\
& =\frac{(aq)_N(aq/b_kc_k)_N}{(aq/b_k)_N(aq/c_k)_N}\sum_{m_1,\ldots, m_{k-1}\geq 0}
\frac{(aq/b_1c_1)_{m_1}(aq/b_2c_2)_{m_2}\cdots (aq/b_{k-1}c_{k-1})_{m_{k-1}}}
{(q)_{m_1}(q)_{m_2}\cdots (q)_{m_{k-1}}}\nonumber\\[3pt]
&\quad\cdot \frac{(b_2)_{m_1}(c_2)_{m_1}(b_3)_{m_1+m_2}(c_3)_{m_1+m_2}\cdots (b_k)_{m_1+\cdots+m_{k-1}}}
{(aq/b_1)_{m_1}(aq/c_1)_{m_1}(aq/b_2)_{m_1+m_2}(aq/c_2)_{m_1+m_2}\cdots (aq/b_{k-1})_{m_1+\cdots+m_{k-1}}}\nonumber\\[3pt]
&\quad\cdot \frac{(c_k)_{m_1+\cdots+m_{k-1}}}{(aq/c_{k-1})_{m_1+\cdots+m_{k-1}}}\cdot
\frac{(q^{-N})_{m_1+m_2+\cdots +m_{k-1}}}{(b_kc_kq^{-N}/a)_{m_1+m_2+\cdots +m_{k-1}}}\nonumber\\[3pt]
&\quad \cdot \frac{(aq)^{m_{k-2}+2m_{k-3}+\cdots +(k-2)m_1}q^{m_1+m_2+\cdots +m_{k-1}}}
{(b_2c_2)^{m_1}(b_3c_3)^{m_1+m_2}\cdots (b_{k-1}c_{k-1})^{m_1+m_2+\cdots +m_{k-2}}}.
\end{align}

The following summation formula
can be derived from the above  transformation formula of Andrews. It
can be considered as a Rogers-Ramanujan type
identity for overpartitions.

\begin{thm}\label{even3} For $k\geq 2$, we have
\begin{align}\label{id1}
&\sum_{N_1 \geq N_2 \geq \cdots \geq N_{k-1}\geq 0}
  \frac{q^{N_1(N_1+1)/2+N_2^2+\cdots+N_{k-1}^2+N_2+\cdots+N_{k-1}}(-q;q)_{N_1}}
{(q;q)_{N_1-N_2}\cdots (q;q)_{N_{k-2}-N_{k-1}}(q;q)_{N_{k-1}}}\nonumber\\
&\quad=\frac{(-q;q)_{\infty}(q;q^{2k})_\infty(q^{2k-1};q^{2k})_{\infty}(q^{2k};q^{2k})_\infty}{(q;q)_\infty}.
\end{align}
\end{thm}

\pf
Applying the above transformation formula of Andrews  by setting
all variables to infinity except for  $c_k$, $a$ and $q$, we get
\begin{align*}
&\sum_{ N_1 \geq \cdots\geq N_{k-1}\geq 0}
\frac{(c_k)_{N_{1}}a^{N_1+\cdots +N_{k-1}}q^{N_1(N_1+1)/2+N_2^2+\cdots+N_{k-1}^2}}
{(q)_{N_1-N_2}\cdots (q)_{N_{k-2}-N_{k-1}}(q)_{N_{k-1}}(-c_k)^{N_{1}}}\nonumber\\
&\quad=\frac{(aq/c_k;q)_{\infty}}{(a,q)_{\infty}} \sum_{n\geq 0} \frac{(1-aq^{2n})
(a,c_k;q)_{n}a^{kn}q^{kn^2}}{(q,aq/c_k;q)_{n}c_k^n}.
\end{align*}
Setting $a=q$ and $c_k=-q$, we find that
\begin{align}
&\sum_{ N_1 \geq \cdots\geq N_{k-1}\geq 0}
\frac{q^{N_1(N_1+1)/2+N_2^2+\cdots+N_{k-1}^2+N_2+\cdots +N_{k-1}}(-q)_{N_{1}}}
{(q)_{N_1-N_2}\cdots (q)_{N_{k-2}-N_{k-1}}(q)_{N_{k-1}}}\nonumber \\
&\quad= \frac{(-q;q)_{\infty}}{(q;q)_{\infty}} \sum_{n\geq 0} (-1)^n(1-q^{2n+1})q^{kn^2+(k-1)n}. \label{a1}
\end{align}
Using Jacobi's triple product identity, we get
\begin{align}
&(q;q^{2k})_{\infty}(q^{2k-1};q^{2k})_{\infty}(q^{2k};q^{2k})_{\infty}\nonumber\\
&\quad =\sum_{n=-\infty}^{\infty}(-1)^nq^{kn^2+(k-1)n}\nonumber\\
&\quad =\sum_{n=0}^{\infty}(-1)^n(1-q^{2n+1})q^{kn^2+(k-1)n}.\label{j1}
\end{align}
In view of (\ref{a1}) and (\ref{j1}), we obtain (\ref{id1}).
 This completes the proof.
\qed

Our second Rogers-Ramanujan type identity for
overpartitions is stated as follows.

\begin{thm}\label{odd2}
For $k\geq 2$, we have
\begin{align}
&\sum_{N_1 \geq N_2 \geq \cdots \geq N_{k-1} \geq 0}
  \frac{q^{N_1(N_1+1)/2+N_2^2+\cdots+N_{k-1}^2+N_2+\cdots+N_{k-1}}(-q;q)_{N_1}}
{(q;q)_{N_1-N_2}\cdots(q;q)_{N_{k-2}-N_{k-1}} (q;q)_{N_{k-1}}(-q;q)_{N_{k-1}}}\nonumber\\
&\quad =\frac{(-q;q)_{\infty}(q;q^{2k-1})_\infty(q^{2k-2};q^{2k-1})_{\infty}(q^{2k-1};q^{2k-1})_\infty}{(q;q)_\infty}.\label{id2}
\end{align}
\end{thm}

\pf
Applying Andrews'   transformation formula by setting  all variables except for
 $c_1$, $c_k$, $a$ and  $q$  to infinity, we find \begin{align*}
&\sum_{ N_1 \geq \cdots\geq N_{k-1}\geq 0}
\frac{(c_k)_{N_{1}}a^{N_1+\cdots +N_{k-1}}q^{N_1(N_1+1)/2+N_2^2+\cdots+N_{k-1}^2}}
{(q)_{N_1-N_2}\cdots (q)_{N_{k-2}-N_{k-1}}(q)_{N_{k-1}}(-c_k)^{N_{1}}(aq/c_1)_{N_{k-1}}}\\
&\quad =\frac{(aq/c_k;q)_{\infty}}{(a,q)_{\infty}}
\sum_{n\geq 0}
\frac{(-1)^n(1-aq^{2n})(a,c_k;q)_{n}(c_1)_n a^{kn}q^{kn^2-(n-1)n/2}}{(q,aq/c_k;q)_{n}(aq/c_1)_n c_1^nc_k^n}.
\end{align*}
Moreover, setting $a=q$, $c_k=-q$ and  $c_1=-q$ yields
\begin{align}
&\sum_{ N_1 \geq \cdots\geq N_{k-1}\geq 0}
\frac{q^{N_2+\cdots +N_{k-1}}q^{N_1(N_1+1)/2+N_2^2+\cdots+N_{k-1}^2}(-q)_{N_{1}}}
{(q)_{N_1-N_2}\cdots (q)_{N_{k-2}-N_{k-1}}(q)_{N_{k-1}}(-q)_{N_{k-1}}}\nonumber\\
&\quad =\frac{(-q;q)_{\infty}}{(q,q)_{\infty}}
\sum_{n\geq 0}
(-1)^n(1-q^{2n+1}) q^{kn^2+kn-n^2/2-3n/2}.\label{a2}
\end{align}
Using Jacobi's triple product identity, we have
\begin{align}
&(q;q^{2k-1})_{\infty}(q^{2k-2};q^{2k-1})_{\infty}(q^{2k-1};q^{2k-1})_{\infty}\nonumber\\
&\quad =\sum_{n=-\infty}^{\infty}(-1)^nq^{kn^2+kn-n^2/2-3n/2}\nonumber\\
&\quad = \sum_{n=0}^{\infty}(-1)^n(1-q^{2n+1})q^{kn^2+kn-n^2/2-3n/2}.\label{j2}
\end{align}
Combining (\ref{a2}) and (\ref{j2}), we deduce (\ref{id2}).
This complete the proof.
\qed

\section{The case when $k$ is even}

In this section, we shall give a proof of
Theorem \ref{main} for an even upper bound $2k-2$. More precisely,  this case can be stated as follows.

\begin{thm}\label{thmeven}For $k\geq 2$ and $n\geq 1$, we have
\begin{equation}
 |F_{2k-2}(n)|=|H_{2k}(n)|.
\end{equation}
\end{thm}

 Recall that the generating function for overpartitions in  $H_{2k}(n)$ equals
 \begin{equation}
 \sum_{n\geq 0 }|H_{2k}(n)|q^n=\frac{(-q;q)_{\infty}(q;q^{2k})_\infty(q^{2k-1};q^{2k})_{\infty}
 (q^{2k};q^{2k})_\infty}{(q;q)_\infty}.
 \end{equation}

 In view of ({\ref{id1}}),  in order to prove Theorem \ref{thmeven} we only need to show that the generating function of anti-lecture hall compositions in
 $F_{2k-2}(n)$
 equals
 the left hand side of ({\ref{id1}}), as stated below.

\begin{thm}\label{even1}
The generating function of anti-lecture hall compositions in  $F_{2k-2}(n)$ is given by
\begin{equation}
 \sum_{n=0}^{\infty}|F_{2k-2}(n)|q^n
=\sum_{N_1 \geq N_2 \geq \cdots \geq N_{k-1}\geq 0}
  \frac{q^{N_1(N_1+1)/2+N_2^2+\cdots+N_{k-1}^2+N_2+\cdots+N_{k-1}}(-q;q)_{N_1}}
{(q;q)_{N_1-N_2}\cdots (q;q)_{N_{k-2}-N_{k-1}}(q;q)_{N_{k-1}}}.
\end{equation}
\end{thm}

In order to prove Theorem \ref{even1}, we  introduce a triangular representation
$T(\lambda)=(t_{ij})_{1 \leqslant i \leqslant j}$
 of an anti-lecture hall composition $\lambda$ which is similar to a T-triangles introduced by
Bousquet-M\'{e}lou \cite{bous99}.

 It should be noted that Corteel and Savage \cite{cor02} used a representation of a composition $\lambda$ as
 a pair of vectors
$(l,r)=((l_1, l_2, \ldots),(r_1,r_2,\ldots))$, where $\lambda_i=il_i+r_i$, with $0\leq r_i\leq i-1$.
Then $l=\lfloor \lambda \rfloor=(\lfloor \lambda_1/1 \rfloor, \lfloor \lambda_2/2 \rfloor,\ldots)$.
It can be checked that a composition $\lambda$ is an anti-lecture hall composition if and only if
\begin{itemize}
\item[(1)] $l_1\geq l_2 \geq \cdots \geq 0$, and

\item[(2)] $r_{i}\geq r_{i+1}$ whenever $l_i=l_{i+1}$.
\end{itemize}

\begin{defi}
The A-triangular representation $T(\lambda)=(t_{i,j})_{1 \leqslant i \leqslant j}$
of an anti-lecture hall composition $\lambda=(\lambda_1,\lambda_2,\ldots)$ is defined to be a
triangular array $(t_{i,j})_{1 \leqslant i \leqslant j}$
of nonnegative integers satisfying the following conditions:
\begin{itemize}
\item[(1)] A diagonal entry $t_{j,j}$ in $T(\lambda)$ equals $l_j=\lfloor \lambda_{j}/j\rfloor$.

\item[(2)]  The first $r_j$ entries of the $j$-th column are equal to $t_{j,j}+1$, while
the other entries in the $j$-th column are equal to $t_{j,j}$.
\end{itemize}
\end{defi}
The sum of all entries of $T(\lambda)$ is equal to $|\lambda|=\lambda_1+\lambda_2+\cdots$. It can be verified that
 the A-triangular representation $T(\lambda)$
of an anti-lecture hall composition possesses the following properties:
\begin{itemize}
\item[(1)]  The diagonal entries of $T$ are weakly decreasing, that is,
 $t_{1,1}\geq t_{2,2} \geq \cdots \geq 0$.

\item[(2)] The entries in the $j$-th column are non-increasing, and they are equal to either the $t_{j,j}$ or $t_{j,j}+1$.

\item[(3)] If $t_{j,j}=t_{j+1,j+1}$, then $t_{i,j}\geq t_{i,j+1}$.
\end{itemize}

Conversely, a triangular
 array satisfying the above conditions must  be the A-triangular representation of an anti-lecture
hall composition.

For example, let
 $\lambda=(4,8,11,14,16,15,11,10,5,2)$. The A-triangular  representation $T(\lambda)$ of $\lambda$  is illustrated
as follows.

\begin{center}
\begin{tabular}{llllllllll}
4&4&4&4&4&3&2&2&1&1\\
&4&4&4&3&3&2&2&1&1\\
&&3&3&3&3&2&1&1&0\\
&&&3&3&2&2&1&1&0\\
&&&&3&2&1&1&1&0\\
&&&&&2&1&1&0&0\\
&&&&&&1&1&0&0\\
&&&&&&&1&0&0\\
&&&&&&&&0&0\\
&&&&&&&&&0\\
\end{tabular}
\end{center}

Now we are ready to give a proof of Theorem \ref{even1} by using the A-triangular representation of an anti-lecture hall
composition.

\noindent{\it Proof of Theorem \ref{even1}.}
Let $\lambda$ be an anti-lecture hall composition  with $\lambda_1\leq 2k-2$.  Let us consider the
 A-triangular representation $T(\lambda)$ of
 $\lambda$.
We use $N_i$ to denote the number of diagonal entries $t_{j,j}$ in $T(\lambda)$ which are greater than or
 equal to
 $2i-1$ for $1\leq i \leq k-1$. Then we have
$N_1\geq N_2 \geq \cdots \geq N_{k-1}\geq 0$. Let $F_{2k-2}(N_1, \ldots, N_{k-1}; n)$ denote
the set of anti-lecture hall compositions $\lambda$ such that there are $N_i$ diagonal entries in $T(\lambda)$
that are greater than or equal to $2i-1$ and $\lambda_1\leq 2k-2$. We aim to compute the generating function of
anti-lecture hall composition in
 $F_{2k-2}(N_1, \ldots, N_{k-1}; n)$, which can be summed up to yield
the generating function of the anti-lecture hall compositions in  $F_{2k-2}(n)$.

Let  $\lambda$ be an anti-lecture hall composition  in $F_{2k-2}(N_1, \ldots, N_{k-1}; n)$, and
let $\lambda^{(1)}=(\lambda_1,\ldots,$ $\lambda_{N_1})$, $\lambda^{(2)}=(\lambda_{N_1+1},\ldots,\lambda_l)$.
Since $\lfloor \lambda_{N_1+1}/(N_1+1)\rfloor=\cdots=\lfloor \lambda_{l}/l\rfloor=0$,   we see that
$\lambda_l\leq \cdots \leq \lambda_{N_1+1} \leq N_1$.
Evidently $\lambda^{(2)}$ is a partition  whose first part is less than $N_1+1$,
and the generating function for possible choices of $\lambda^{(2)}$ equals $1/(q;q)_{N_1}$.

 Let us examine the composition $\lambda^{(1)}$ and its A-triangular representation $T(\lambda^{(1)})$. The triangular array $T(\lambda^{(1)})$
can be split into  $k$ triangular arrays and we can  compute the
generating function for possible choices of $\lambda^{(1)}$.

\noindent
Step 1. Let $T^{(1)}=T(\lambda^{(1)})$. Extract $1$ from each entry in the first $N_1$ columns of $T^{(1)}$ to form a
triangular array of size $N_1$ with all the
entries equal to $1$, denoted by $R(N_1,1)$.

\noindent
Step 2. For $2\leq i \leq k-1$, extract $2$ from each entry in the first $N_i$ columns of the remaining
triangular array $T^{(1)}$  to
     generate a triangular array of size $N_i$ with all the entries equal to $2$, denoted by $R(N_i,2)$.

\noindent
Step 3. Let $S$ denote the remaining triangular array $T^{(1)}$.

After the above operations, $T(\lambda^{(1)})$ is decomposed  into $k$ triangular arrays, including an A-triangle $R(N_1,1)$
of size $N_1$ with  entries $1$,
$k-2$ A-triangular arrays $R(N_i,2)$ of sizes $N_2, \ldots, N_{k-1}$ respectively
 with  entries $2$ where  $i=2, \ldots, k-1$, and
a triangular array $S=(s_{i,j})_{1\leq i \leq j\leq N_1}$ of  size $N_1$. It is easy to see that the generating function
 for triangular arrays in $R(N_1,1)$ is
$q^{(N_1+1)N_1/2}$ and the generating function
of triangular arrays in $R(N_i,2)$ is $q^{N_i^2+N_i}$.

It can be verified that $S$ possesses the following properties by the definition of the A-triangular representation of an
anti-lecture hall composition:
\begin{itemize}
\item[(1)] All the entries in the diagonals of $S$ are equal  to $1$ or $0$. Note that $S$ has $N_1$ diagonal
elements
$s_{1,1}, s_{2,2}, \ldots, s_{N_1,N_1}$.
    These diagonal elements can be divided into $k-1$ segments such that  the first
     segment contains $n_1=N_1-N_2$ elements $s_{N_2+1,N_2+1},\ldots, s_{N_1,N_1}$,
the second segment contains $n_2=N_2- N_3$ elements $s_{N_3+1,N_3+1},\ldots, s_{N_2,N_2}$,
 and so on, while the last segment contains $n_{k-1}=N_{k-1}$ elements $s_{1,1},\ldots,s_{N_{k-1},N_{k-1}}$.
Moreover, the $i$-th segment contains $m_i$ 1's followed by 0's.

\item[(2)] The entries in the $j$-th column are non-increasing, and they are equal to
either the $t_{j,j}$ or $t_{j,j}+1$.

\item[(3)] If $s_{j,j}=s_{j+1,j+1}$, then $s_{i,j}\geq s_{i,j+1}$.
\end{itemize}

We denote the set of triangular arrays possessing the above three properties by $S(N_1, N_2, \ldots,$ $N_{k-1})$. Now we are in a position to compute
 the generating function of triangular arrays in $S(N_1,N_2,\ldots,N_{k-1})$.

 We may partition a triangular array $S\in S(N_1,N_2,\ldots,N_{k-1})$   into $k-1$ blocks of columns, where
  the $i$-th block consists of the
 $(N_{i+1}+1)$-th column to the $N_i$-th column of $S$. We denote the $i$-th block by $S_i$.
According to the above three properties, we deduce that the first $m_i$ diagonal entries of $S_i$ must be $1$ and the entries in the first $m_i$ columns of $S_i$
  are either $1$ or $2$.

We shall split $S_i$ into three trapezoidal arrays $S_i^{(1)}$, $S_i^{(2)}$ and $S_i^{(3)}$.
First, we may form a  trapezoidal array $S_i^{(1)}$  of the same size as $S_i$ and with the entries
in the first $m_i$ columns  equal  to $1$ and the other entries equal to 0.
     Let $S^{\,'}_i$ denote the trapezoidal  array obtained from $S_i$ by subtracting
     1 from every entry in the first $m_i$ columns. Observe that every entry in
      $S^{\,'}_i$ is either $1$ or $0$, and $S^{(1)}_i$ can be regarded as the Ferrers diagram
 of the conjugate of the partition
\[ \alpha^{(1)}=(N_{i+1}+m_i, N_{i+1}+m_i-1, \ldots, N_{i+1}+1).\]
Furthermore,   $S^{\,'}_i$ satisfies the following conditions:
\begin{itemize}
\item[(1)] All entries in $S^{\,'}_i$ are equal to $0$ or $1$, but the diagonal entries   must be $0$.

\item[(2)] The entries in the same column must be non-increasing.

\item[(3)] The first $m_i$ entries in the $j$-th row must be non-increasing,  and the remaining entries in the $j$-th
 row are also non-increasing.
\end{itemize}

We continue to consider the trapezoidal array
formed by the first $m_i$ columns of  $S^{\,'}_i$, and
denote it by $S_i^{(2)}$.
Similarly, we see that $S_i^{(2)}$ can be
 regarded
as the Ferrers diagram of the conjugate of a partition $\alpha^{(2)}$, where \[ \alpha^{(2)}_1\leq N_{i+1}, \quad \mbox{and} \quad
 l(\alpha^{(2)})\leq m_i.\]

 Define  $S_i^{(3)}$ to be the trapezoidal array
formed by the
 $(m_i+1)$-th column to the $(N_i-N_{i+1})$-th column  of $S'_i$. Again, $S^{(3)}$ can be regarded
as the Ferrers diagram of the conjugate of a partition $\alpha^{(3)}$, where
 \[ \alpha^{(3)}_1\leq N_{i+1}+m_i\quad \mbox{and} \quad
l(\alpha^{(3)})\leq N_i-N_{i+1}-m_i.\] So the generating function for  possible choices of the $i$-th block $S_i$ is given by
\begin{align}\label{si1}
\sum_{m_i=0}^{N_i-N_{i+1}}q^{\frac{(N_{i+1}+1+N_{i+1}+m_i)m_i}{2}}
\frac{(q;q)_{N_{i+1}+m_i}}{(q;q)_{m_i}(q;q)_{N_{i+1}}}
\frac{(q;q)_{N_i}}{(q;q)_{N_{i+1}+m_i}(q;q)_{N_i-N_{i+1}-m_i}}.
\end{align}
which equals
\begin{align}
\frac{(q;q)_{N_i}}{(q;q)_{N_{i+1}}(q;q)_{N_i-N_{i+1}}}\sum_{m_i=0}^{N_i-N_{i+1}}
q^{\frac{(N_{i+1}+1+N_{i+1}+m_i)m_i}{2}}
\frac{(q;q)_{N_i-N_{i+1}}}{(q;q)_{m_i}(q;q)_{N_i-N_{i+1}-m_i}}.
\end{align}
Observe that the sum
\[\sum_{m_i=0}^{N_i-N_{i+1}}q^{\frac{(N_{i+1}+1+N_{i+1}+m_i)m_i}{2}}
\frac{(q;q)_{N_i-N_{i+1}}}{(q;q)_{m_i}(q;q)_{N_i-N_{i+1}-m_i}}\]
 is   the generating function for    partitions
with distinct parts between $N_{i+1}+1$ and $N_i$. Therefore,
\begin{align}\label{si2}
 \sum_{m_i=0}^{N_i-N_{i+1}}q^{\frac{(N_{i+1}+1+N_{i+1}+m_i)m_i}{2}}
\frac{(q;q)_{N_i-N_{i+1}}}{(q;q)_{m_i}(q;q)_{N_i-N_{i+1}-m_i}}
=(-q^{N_{i+1}+1};q)_{N_i-N_{i+1}}.
\end{align}
By (\ref{si2}),  the generating function (\ref{si1}) can be simplified to
\begin{equation}
\frac{(q;q)_{N_i}}{(q;q)_{N_{i+1}}(q;q)_{N_i-N_{i+1}}}(-q^{N_{i+1}+1};q)_{N_i-N_{i+1}}.
\end{equation}
Thus the generating function for triangular arrays in $S$  can be written as
\begin{equation*}
\prod_{i=1}^{k-1}\frac{(q;q)_{N_i}}{(q;q)_{N_{i+1}}(q;q)_{N_i-N_{i+1}}}(-q^{N_{i+1}+1};q)_{N_i-N_{i+1}}
=\frac{(q)_{N_1}(-q;q)_{N_1}}{(q)_{N_1-N_2}\cdots (q)_{N_{k-2}-N_{k-1}}(q)_{N_{k-1}}}.
\end{equation*}

Recall that the generating function for  possible choices of  $T(\lambda^{(2)})$ equals $1/(q;q)_{N_1}$ and
the generating functions for $R(N_1,1)$, $R(N_2,2), \ldots, R(N_{k-1},2)$ are equal to $q^{(N_1+1)N_1/2}$, $q^{N_2^2+N_2}$,$\ldots,$ $q^{N_{k-1}^2+N_{k-1}}$
respectively. We also note that
the generating function for anti-lecture hall compositions in
 $F_{2k-2}(N_1,\ldots, N_{k-1}, n)$ is  the product of the generating functions for
$T(\lambda^{(2)})$, $R(N_1,1)$, $R(N_2,2), \ldots, R(N_{k-1},2)$
and $S$, and therefore it equals
\begin{align*}
&\frac{q^{(N_1+1)N_1/2+N^2_2+\cdots+N_{k-1}^2+N_2+\cdots+N_{k-1}}}{(q)_{N_1}}
\frac{(q)_{N_1}(-q;q)_{N_1}}{(q)_{N_1-N_2}\cdots (q)_{N_{k-2}-N_{k-1}}(q)_{N_{k-1}}}\\
&\quad=
\frac{q^{(N_1+1)N_1/2+N^2_2+\cdots+N_{k-1}^2+N_2+\cdots+N_{k-1}}(-q;q)_{N_1}}
{(q)_{N_1-N_2}\cdots (q)_{N_{k-2}-N_{k-1}}(q)_{N_{k-1}}}.
\end{align*}

Summing up the generating functions of anti-lecture hall compositions in
 $F_{2k-2}(N_1,\ldots, N_{k-1}, n)$, we get the generating function for anti-lecture hall compositions in  $F_{2k-2}(n)$,
\begin{equation}
\sum_{n\geq 0}|F_{2k-2}(n)|q^n
=\sum_{N_1\geq \cdots \geq N_{k-1}\geq 0}
\frac{q^{(N_1+1)N_1/2+N^2_2+\cdots+N_{k-1}^2+N_2+\cdots+N_{k-1}}(-q;q)_{N_1}}
{(q)_{N_1-N_2}\cdots (q)_{N_{k-2}-N_{k-1}}(q)_{N_{k-1}}}.
\end{equation}
The proof is therefore completed.
\qed

For example, let  $\lambda=(4,8,11,14,16,15,11,10,5,2)$ and let $k=3$. Then $\lambda^{(2)}=(5,2)$, $N_1=8$, $N_2=5$,
$m_1=2$ and $m_2=1$.
The decomposition of $T(\lambda)$
is illustrated as follows:
\begin{center}
\setlength{\tabcolsep}{5.5pt}
\begin{tabular}{ccccc}
\begin{tabular}{llllllllll}
4&4&4&4&4&3&2&2&1&1                \\
&4&4&4&3&3&2&2&1&1                     \\
&&3&3&3&3&2&1&1&0           \\
&&&3&3&2&2&1&1&0            \\
&&&&3&2&1&1&1&0       \\
&&&&&2&1&1&0&0            \\
&&&&&&1&1&0&0             \\
&&&&&&&1&0&0            \\
&&&&&&&&0&0               \\
&&&&&&&&&0
\end{tabular}
&\hspace{-10pt}$\longrightarrow$
&
\hspace{-10pt}\begin{tabular}{llllllll}
         1&1&1&1&1&1&1&1                    \\
         &1&1&1&1&1&1&1             \\
         &&1&1&1&1&1&1             \\
         &&&1&1&1&1&1            \\
         &&&&1&1&1&1  \\
         &&&&&1&1&1                 \\
         &&&&&&1&1                  \\
         &&&&&&&1                 \\
         &&&&&&&                   \\
         &&&&&&                     \\
\end{tabular}
&
\hspace{-10pt}$+$
&
\hspace{-10pt}\begin{tabular}{lllll}
       2&2&2&2&2           \\
       &2&2&2&2\\
       &&2&2&2\\
       &&&2&2\\
       &&&&2\\
       &&&& \\
        &&&& \\
 &&&& \\
 &&&& \\
 &&&&
 \end{tabular}\\

 $T(\lambda)$& & $R(8,1)$ & & $R(5,2)$\\
\end{tabular}
\end{center}

\setlength{\tabcolsep}{5.5pt}

\vspace{10pt}\hspace{40pt}\begin{tabular}{cccc}
$+$
&
\begin{tabular}{lllllllll}
   &1&1&1&1&1&2&1&1                          \\
   &&1&1&1&0&2&1&1                     \\
   &&&0&0&0&2&1&0                      \\
   &&&&0&0&1&1&0                           \\
&&&&&0&1&0&0            \\
   &&&&&&1&0&0                             \\
   &&&&&&&0&0                                \\
   &&&&&&&&0                                  \\
   &&&&&&&&         \\
   &&&&&&&&
\end{tabular}
&
$+$
&
\begin{tabular}{llllllllll}
   0&0&0&0&0&0&0&0&1&1                             \\
   &0&0&0&0&0&0&0&1&1                     \\
   &&0&0&0&0&0&0&1&0                    \\
   &&&0&0&0&0&0&1&0                          \\
   &&&&0&0&0&0&1&0           \\
   &&&&&0&0&0&0&0                            \\
   &&&&&&0&0&0&0                               \\
   &&&&&&&0&0&0                                  \\
   &&&&&&&&0&0                                    \\
   &&&&&&&&&0
\end{tabular}\\

&$S$&&$T(\lambda^{(2)})$

\end{tabular}

\section{The case when $k$ is odd}

 The objective of this section is to provide two proofs of the following theorem which is
the  case of an odd upper bound $2k-3$ of Theorem \ref{main}. The first
is analogous to the proof of the even case.
The second  requires a Rogers-Ramanujan type identity of Andrews, a bijection of Corteel and Savage, and a refined version of a bijection also due to  Corteel and
Savage. The approach of the second proof does not seem to apply to the even case, namely, Theorem \ref{thmeven}.

\begin{thm}\label{thmodd}
For $k\geq 2$ and a positive integer $n$, we have
\begin{equation}
|F_{2k-3}(n)|=|H_{2k-1}(n)|.
\end{equation}
\end{thm}

The first proof relies on the following
  generating function formula for anti-lecture hall compositions in $F_{2k-3}(n)$.
  The proof of this formula is analogous to that of Theorem \ref{even1}.

\begin{thm}\label{odd1}
For $k\geq 2$,
\begin{equation}\label{gaodd}
\sum_{n=0}^{\infty}|F_{2k-3}(n)|q^n=\sum_{N_1 \geq N_2 \geq \cdots \geq N_{k-1} \geq 0}
  \frac{q^{(N_1+1)N_1/2+N_2^2+\cdots+N_{k-1}^2+N_2+\cdots+N_{k-1}}(-q;q)_{N_1}}
{(q;q)_{N_1-N_2}\cdots(q;q)_{N_{k-2}-N_{k-1}} (q;q)_{N_{k-1}}(-q;q)_{N_{k-1}}}.
\end{equation}
\end{thm}

\noindent
{\it Proof of Theorem \ref{odd1}.}
 Let $\lambda$ be an anti-lecture hall composition  with $\lambda_1\leq 2k-3$.  We consider the
 A-Triangular representation $T(\lambda)$ of
 $\lambda$.
Let $N_i$ be the number of diagonal entries $t_{jj}$ in $T(\lambda)$ which are greater than or equal to $2i-1$ for
$1\leq i \leq k-1$. Then we have
$N_1\geq N_2 \geq \cdots \geq N_{k-1}\geq 0$. Let $F_{2k-3}(N_1, \ldots, N_{k-1}; n)$ denote
the set of anti-lecture hall compositions $\lambda$ for which there are $N_i$ diagonal entries in $T(\lambda)$
that are greater than or equal to  $2i-1$ and $\lambda_1\leq 2k-3$.

Let $\lambda^{(1)}=(\lambda_1,\ldots,\lambda_{N_1})$, $\lambda^{(2)}=(\lambda_{N_1+1},\ldots,\lambda_l)$.
It is immediately verified that $\lambda^{(2)}$ is a partition  whose first part does not exceed $N_1$. Hence the
 generating function for possible choices of
$\lambda^{(2)}$ equals $1/(q;q)_{N_1}$.

Now consider $\lambda^{(1)}$ and its A-Triangular representation $T(\lambda^{(1)})$. We can split $T(\lambda^{(1)})$
into $k$ triangular arrays to  compute the
generating function for possible choices of  $\lambda^{(1)}$.

\noindent
Step 1. Let $T^{(1)}=T(\lambda^{(1)})$. Extract $1$ from each entry in the first $N_1$ columns of $T^{(1)}$ to form a
 triangular array of size $N_1$ with all
entries equal to $1$, denoted by $R(N_1,1)$.

\noindent
Step 2. For $i=2, \ldots, k-1$, extract $2$ from each entry in the first $N_i$ columns of the remaining array $T^{(1)}$  to
      form a triangular array of size $N_i$ with all entries equal to  $2$, denoted by $R(N_i,2)$.

\noindent
Step 3. Let $S$ be the remaining triangular array $T^{(1)}$.

After the above procedures, $T(\lambda^{(1)})$ is decomposed into $k$ triangular arrays, including an A-Triangle $R(N_1,1)$
of size $N_1$ with all entries being $1$,
$(k-2)$ A-Triangles $R(N_i,2)$ of sizes $N_2,\ldots, N_{k-1}$ respectively with  all entries being $2$ and
a triangular array $S=(s_{i,j})$ of size $N_1$
satisfying the following conditions:
\begin{itemize}

  \item[(1)]
All the entries in the diagonals of $S$ are equal  to $1$ or $0$. Note that $S$ has $N_1$ diagonal elements
$s_{1,1}, s_{2,2}, \ldots, s_{N_1,N_1}$.
    These diagonal elements can be divided into $k-1$ segments such that  the first
     segment contains $n_1=N_1-N_2$ elements $s_{N_2+1,N_2+1},\ldots, s_{N_1,N_1}$,
the second segment contains $n_2=N_2- N_3$ elements $s_{N_3+1,N_3+1},\ldots, s_{N_2,N_2}$,
 and so on, while the last segment contains $n_{k-1}=N_{k-1}$ elements $s_{1,1},\ldots,s_{N_{k-1},N_{k-1}}$.
Moreover, the $i$-th segment contains $m_i$ 1's followed by 0's.

\item[(2)] The entries in the $j$-th column are non-increasing, and they are equal to
either  $t_{j,j}$ or $t_{j,j}+1$.

\item[(3)] If $s_{j,j}=s_{j+1,j+1}$, then $s_{i,j}\geq s_{i,j+1}$.

\item[(4)] The entries in the first $N_{k-1}$ columns of $S$ are equal to $0$, that is, $m_{k-1}=0$.
\end{itemize}

Let us write $\overline{S}(N_1,N_2,\cdots, N_{k-1})$ for the set of triangular arrays possessing the above four properties.
We proceed to compute the generating function for the triangular arrays in $\overline{S}(N_1,N_2,\cdots$ $,N_{k-1})$.

We may partition a triangular array $S\in \overline{S}(N_1,N_2,\ldots,N_{k-1})$   into $k-1$ blocks of columns, where
  the $i$-th block consists of the
 $(N_{i+1}+1)$-th column to the $N_i$-th column of $S$. We denote the $i$-th block by $S_i$.
According to the above four properties, we infer that the first $m_i$ diagonal entries of $S_i$ must be $1$ and
the entries in the first $m_i$ columns of $S_i$
  are either $1$ or $2$ for $i=1,\ldots k-2$ and $S_{k-1}$ is a triangular array of size $N_{k-1}$ with all entries
eqaul to  $0$.

 We shall split $S_i$ into three trapezoidal arrays $S_i^{(1)}$, $S_i^{(2)}$ and $S_i^{(3)}$ for $i=1,\ldots k-2$.
First, we  may form a  trapezoidal array $S_i^{(1)}$ of the same size as $S_i$ and with the entries
in the first $m_i$ columns  equal  to $1$ and the other entries equal to 0.
     Let $S^{\,'}_i$ denote the trapezoidal  array obtained from $S_i$ by subtracting
     1 from every entry in the first $m_i$ columns. It is seen  that every entry in
      $S^{\,'}_i$ is either $1$ or $0$, and $S^{(1)}_i$ can be regarded as the Ferrers diagram
 of the conjugate of the partition
\[ \alpha^{(1)}=(N_{i+1}+m_i, N_{i+1}+m_i-1, \ldots, N_{i+1}+1).\]

Furthermore,   $S^{\,'}_i$ satisfies the following conditions for $i=1,\ldots, k-2$:
\begin{itemize}
\item[(1)] All the entries in $S^{\,'}_i$ equal $0$ or $1$, but the diagonal entries must be $0$.

\item[(2)] The entries in the $j$-th column must be non-increasing.

\item[(3)] The first $m_i$ entries in the $j$-th row must be non-increasing,  and the remaining entries in the $j$-th
 row are also non-increasing.
\end{itemize}

We continue to consider the trapezoidal array
formed by the first $m_i$ columns of  $S^{\,'}_i$, and
denote it by $S_i^{(2)}$.
Again, we see that $S_i^{(2)}$ can be
 regarded
as the Ferrers diagram of the conjugate of a partition $\alpha^{(2)}$, where \[ \alpha^{(2)}_1\leq N_{i+1}, \quad \mbox{and} \quad
 l(\alpha^{(2)})\leq m_i.\]

 Notice that there are still some columns
 to be dealt with. Define  $S_i^{(3)}$ to be the trapezoidal array
formed by the
 $(m_i+1)$-th column to the $(N_i-N_{i+1})$-th column  of $S^{\,'}_i$. Once more, $S_i^{(3)}$ can be regarded
as the Ferrers diagram of the conjugate of a partition $\alpha^{(3)}$, where
 \[ \alpha^{(3)}_1\leq N_{i+1}+m_i\quad \mbox{and} \quad
l(\alpha^{(3)})\leq N_i-N_{i+1}-m_i.\]

 As a consequence, the generating function for possible choices of the $i$-th block $S_i$ for $i=1,\ldots,k-2$ equals
\begin{align*}
\sum_{m_i=0}^{N_i-N_{i+1}}q^{\frac{(N_{i+1}+1+N_{i+1}+m_i)m_i}{2}}\frac{(q;q)_{N_{i+1}+m_i}}{(q;q)_{m_i}(q;q)_{N_{i+1}}}
\frac{(q;q)_{N_i}}{(q;q)_{N_{i+1}+m_i}(q;q)_{N_i-N_{i+1}-m_i}}
\end{align*}
which can be rewritten as
\begin{align*}
\frac{(q;q)_{N_i}}{(q;q)_{N_{i+1}}(q;q)_{N_i-N_{i+1}}}\sum_{m_i=0}^{N_i-N_{i+1}}
q^{\frac{(N_{i+1}+1+N_{i+1}+m_i)m_i}{2}}
\frac{(q;q)_{N_i-N_{i+1}}}{(q;q)_{m_i}(q;q)_{N_i-N_{i+1}-m_i}}.
\end{align*}
Evidently, the sum in the above expression
is  the generating function for  partitions
with distinct parts between $N_{i+1}+1$ and $N_i$. So we deduce that
\begin{align*}
\sum_{m_i=0}^{N_i-N_{i+1}}q^{\frac{(N_{i+1}+1+N_{i+1}+m_i)m_i}{2}}
\frac{(q;q)_{N_i-N_{i+1}}}{(q;q)_{m_i}(q;q)_{N_i-N_{i+1}-m_i}}
=(-q^{N_{i+1}+1};q)_{N_i-N_{i+1}}.
\end{align*}
Since the generating function for $S_{k-1}$ equals $1$, the generating function for possible choices of $S$ is the product of the
generating functions for $S_i$ for $i=1,\ldots,k-2$, that is,
\begin{equation*}
\prod_{i=1}^{k-2}\frac{(q;q)_{N_i}}{(q;q)_{N_{i+1}}(q;q)_{N_i-N_{i+1}}}(-q^{N_{i+1}+1};q)_{N_i-N_{i+1}}
=\frac{(q)_{N_1}(-q;q)_{N_1}}{(q)_{N_1-N_2}\cdots (q)_{N_{k-2}-N_{k-1}}(q)_{N_{k-1}}(-q;q)_{N_{k-1}}}.
\end{equation*}

Recall that the generating function for possible choices of  $T(\lambda^{(2)})$ equals $1/(q;q)_{N_1}$ and
the generating functions
 for  $R(N_1,1)$, $R(N_2,2), \ldots, R(N_{k-1},2)$ are equal to
$q^{(N_1+1)N_1/2}$, $q^{N_2^2+N_2},\ldots$ $, q^{N_{k-1}^2+N_{k-1}}$
respectively. We also observe that
the generating function for anti-lecture hall compositions in  $F_{2k-2}(N_1,\ldots, N_{k-1}, n)$ is  the
 product of the generating functions for
$T(\lambda^{(2)})$, $R(N_1,1)$, $R(N_2,2), \ldots, R(N_{k-1},2)$
and $S$. Hence it equals
\begin{align*}
&\frac{q^{(N_1+1)N_1/2+N^2_2+\cdots+N_{k-1}^2+N_2+\cdots+N_{k-1}}}{(q)_{N_1}}
\frac{(q)_{N_1}(-q;q)_{N_1}}{(q)_{N_1-N_2}\cdots (q)_{N_{k-2}-N_{k-1}}(q)_{N_{k-1}}(-q;q)_{N_{k-1}}}\\
&=\quad
\frac{q^{(N_1+1)N_1/2+N^2_2+\cdots+N_{k-1}^2+N_2+\cdots+N_{k-1}}(-q;q)_{N_1}}
{(q)_{N_1-N_2}\cdots (q)_{N_{k-2}-N_{k-1}}(q)_{N_{k-1}}(-q;q)_{N_{k-1}}}.
\end{align*}
Summing up the generating functions for anti-lecture hall compositions in $F_{2k-3}(N_1,\ldots, N_{k-1}, n)$ yields the generating function for $F_{2k-3}(n)$,
\begin{equation}
 \sum_{n\geq 0}|F_{2k-3}(n)|q^n
=\quad \sum_{N_1\geq \cdots \geq N_{k-1}\geq 0}
\frac{q^{(N_1+1)N_1/2+N^2_2+\cdots+N_{k-1}^2+N_2+\cdots+N_{k-1}}(-q;q)_{N_1}}
{(q)_{N_1-N_2}\cdots (q)_{N_{k-2}-N_{k-1}}(q)_{N_{k-1}}(-q;q)_{N_{k-1}}}.
\end{equation}
This completes the proof.
\qed

For example, the composition $\lambda=(5,10,14,17,18,20,18,15,12,3)$ can decomposed into the following triangular arrays

\begin{center}
\setlength{\tabcolsep}{4.5pt}
\begin{tabular}{cccccccc}

\hspace{-10pt}\begin{tabular}{llllllllll}
5&5&5&5&4&4&3&2&2&1   \\
&5&5&4&4&4&3&2&2&1  \\
&&4&4&4&3&3&2&2&1       \\
&&&4&3&3&3&2&1&0         \\
&&&&3&3&2&2&1&0    \\
&&&&&3&2&2&1&0   \\
&&&&&&2&2&1&0   \\
&&&&&&&1&1&0\\
&&&&&&&&1&0 \\
&&&&&&&&&0
\end{tabular}
&\hspace{-10pt}$\longrightarrow$
&
\hspace{-10pt}\begin{tabular}{lllllllll}
         1&1&1&1&1&1&1&1&1                    \\
         &1&1&1&1&1&1&1&1             \\
         &&1&1&1&1&1&1&1             \\
         &&&1&1&1&1&1&1            \\
         &&&&1&1&1&1&1  \\
         &&&&&1&1&1&1                 \\
         &&&&&&1&1&1                  \\
         &&&&&&&1&1                 \\
         &&&&&&&&1                   \\
         &&&&&&&&                     \\
\end{tabular}
&
\hspace{-10pt}$+$
&
\hspace{-10pt}\begin{tabular}{llllll}
       2&2&2&2&2&2           \\
       &2&2&2&2&2\\
       &&2&2&2&2\\
       &&&2&2&2\\
       &&&&2&2\\
       &&&&&2 \\
        &&&&& \\
 &&&& &\\
 &&&&& \\
 &&&& &
 \end{tabular}
&
\hspace{-10pt}$+$
&
\hspace{-10pt}\begin{tabular}{llllll}
       2&2           \\
       &2\\
       &\\
       &\\
       &\\
       & \\
        & \\
 &\\
 & \\
 &
 \end{tabular}\\

 $T(\lambda)$& & $R(9,1)$ & & $R(6,2)$&&$R(3,2)$\\
\end{tabular}
\end{center}

\setlength{\tabcolsep}{4.5pt}
\vspace{10pt}\hspace{50pt}\begin{tabular}{cccc}
$+$
&
\begin{tabular}{lllllllll}
  0&0&2&2&1&1&2&1&1   \\
&0&2&1&1&1&2&1&1  \\
&&1&1&1&0&2&1&1       \\
&&&1&0&0&2&1&0         \\
&&&&0&0&1&1&0    \\
&&&&&0&1&1&0   \\
&&&&&&1&1&0   \\
&&&&&&&0&0\\
&&&&&&&&0
\end{tabular}
&
$+$
&
 \begin{tabular}{llllllllll}
0&0&0&0&0&0&0&0&0&1   \\
&0&0&0&0&0&0&0&0&1  \\
&&0&0&0&0&0&0&0&1       \\
&&&0&0&0&0&0&0&0         \\
&&&&0&0&0&0&0&0    \\
&&&&&0&0&0&0&0   \\
&&&&&&0&0&0&0   \\
&&&&&&&0&0&0\\
&&&&&&&&0&0 \\
&&&&&&&&&0
\end{tabular}\\

&$S$&&$\lambda^{(2)}$

\end{tabular}

  In virtue of (\ref{id2}),   Theorem \ref{thmodd} immediately follows from Theorem \ref{odd1}, since
   the generating function for overpartitions in
 $H_{2k-1}(n)$ is given by
\begin{align}
\frac{(-q;q)_{\infty}
(q;q^{2k-1})_\infty(q^{2k-2};
q^{2k-1})_{\infty}
(q^{2k-1};q^{2k-1})_\infty}{(q;q)_\infty}.
\end{align}

We now come to the second proof of Theorem \ref{thmodd}. In their proof of anti-lecture hall theorem,
 Corteel and Savage \cite{cor02} established two
 bijections.
The first  is a bijection between
 the set $E(n)$ of  anti-lecture hall compositions $\mu$ of $n$  such that
 $\lfloor \mu_i/i \rfloor$
is even
 and  the set $P(n)$ of
partitions of $n$ with each part greater than one. The  second
bijection is between the set $A(n)$ of anti-lecture hall compositions of $n$  and the set $D \times E(n)$ of pairs $(\lambda, \mu)$ such that $|\lambda|+|\mu|=n$ and
 $\lambda\in D$, $\mu \in E$, where $D$ is the set of partitions into distinct parts.
 Then the anti-lecture hall theorem can
 follows from the correspondence between $A(n)$ and $D \times P(n)$.

 We shall present a  bijection between
a subset of $P(n)$ and a subset of  $E(n)$. Together with the second bijection of
 Corteel and Savage, we arrive at the assertion in Theorem \ref{thmodd}.

To be more specific,
 let $Q_k(n)$ be the subset of $E(n)$ consisting of  anti-lecture hall compositions $\lambda$ such that
$\lambda_1\leq k$ and let $R_k(n)$ be the subset of $P(n)$ consisting of   partitions having at most
$k-1$ successive $N\times (N+1)$ Durfee rectangles such that there is no part below the last Durfee rectangle. Then we have the following
correspondence, which can be considered as a refined version of the first bijection of Corteel and Savage.

\begin{thm}\label{odd3} There is a bijection
between the set
$R_k(n)$ and the set $Q_{2k-2}(n)$.
\end{thm}

\pf
We proceed to give a construction of the bijection $\theta$ from $R_k(n)$ to  $Q_{2k-2}(n)$.
Consider the A-triangular representation $T(\mu)$ of an anti-lecture hall composition $\mu$ of $n$  such that   $\lfloor \frac{\mu_i}{i}\rfloor$ are even for all $i$
and  $\mu_1\leq 2k-2$.  By definition,   all the diagonal entries of $T(\mu)$ are  even and $t_{1,1}\leq 2k-2$.

Now we define the map  $\theta$ from a partition $\lambda$ in $P$
 with exactly $k-1$ successive Durfee rectangles to an
anti-lecture hall composition $\mu$ of $n$.

\noindent
Step 1. We break the Ferrers diagram of $\lambda$ into $k-1$ blocks such that the $i$-th block
contains the $i$-th Durfee rectangle and the dots on the right of  the $i$-th Durfee rectangle.

\noindent
Step 2. Change the $i$-th Durfee rectangle in the $i$-th block into a triangular array with all entries being $2$,
 and the
       rest dots in the $i$-th block into entries equal to $1$. Then these $k-1$ blocks become $k-1$
A-triangles with
all the diagonal entries  equal to $0$ or $2$ where $0$'s are omitted.

\noindent
Step 3. Put  the $k-1$ A-triangles obtained in Step 2 together to form an A-triangle $T$.

The resulting A-triangle
corresponds to an anti-lecture hall composition
$\mu$ such that  $\mu_1=2k-2$ and
 $\lfloor \lambda_i/i \rfloor$ are even for all $i$.

 It is easily verified that the map $\theta$ is reversible. This completes the proof.
\qed

For example,
 let
 \[\lambda=(10,10,9,8,7,7,7,7,5,4,3)\] be a partition in $R_4(77)$. Then the corresponding anti-lecture hall
composition in $Q_6(77)$ equals \[\mu=(6,12,13,11,12,14,4,3,2).\]
 The successive Durfee rectangles of $\lambda$ are exhibited as follows.

\hspace{-30pt} \begin{tabular}{lclcl}
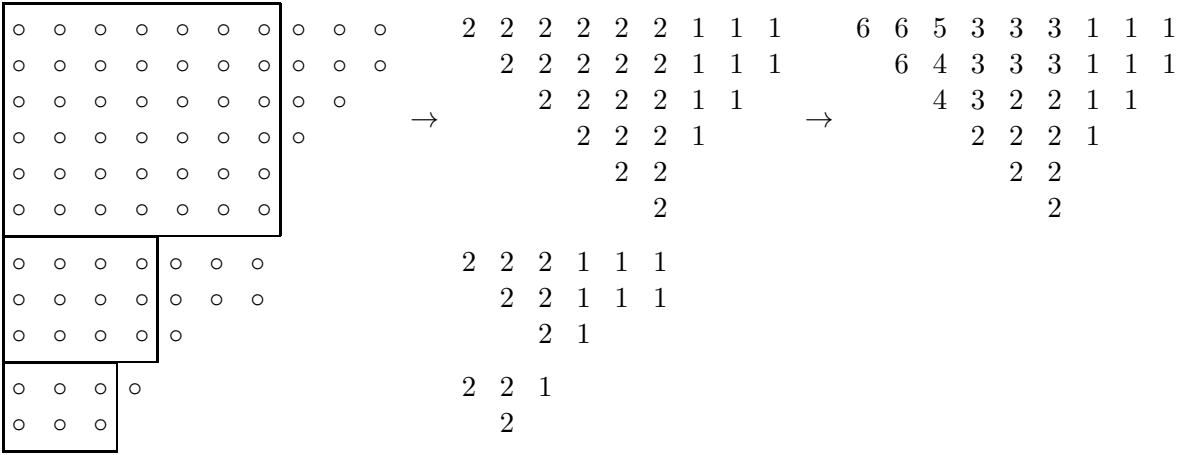
\begin{tabular}{ll}
 \fbox{$\begin{matrix}
 \circ &\circ &\circ &\circ &\circ &\circ &\circ\\
 \circ &\circ &\circ &\circ &\circ &\circ &\circ \\
 \circ &\circ &\circ &\circ &\circ &\circ &\circ \\
 \circ &\circ &\circ &\circ &\circ &\circ &\circ \\
 \circ &\circ &\circ &\circ &\circ &\circ &\circ  \\
 \circ &\circ &\circ &\circ &\circ &\circ &\circ
 \end{matrix}$}
&\hspace{-5pt}$\begin{matrix}
\circ  &\circ &\circ\\
\circ &\circ &\circ\\
\circ &\circ &\\
\circ &&\\
\ &&\\
\ &&\\
  \end{matrix}$
\end{tabular}
&\hspace{-5pt}$\rightarrow$
&\hspace{-5pt}\begin{tabular}{lllllllll}
2 &2 &2 &2 &2 &2&1&1&1\\
 &2 &2 &2 &2 &2&1&1&1\\
&&2 &2 &2 &2&1&1&\\
  && &2 &2 &2&1&&\\
&&&&2&2&&&\\
  &&&&&2&&&
 \end{tabular}
&\hspace{-5pt}$\rightarrow$
&
\hspace{-5pt}\begin{tabular}{lllllllll}
6 &6 &5 &3 &3 &3&1&1&1\\
 &6 &4 &3 &3 &3&1&1&1\\
&&4 &3 &2 &2&1&1&\\
  && &2 &2 &2&1&&\\
&&&&2&2&&&\\
  &&&&&2&&&
 \end{tabular}
\\
\begin{tabular}{ll}
\fbox{$\begin{matrix}
 \circ &\circ &\circ &\circ\\
 \circ &\circ &\circ &\circ\\
 \circ &\circ &\circ &\circ\\
 \end{matrix}$}
&\hspace{-5pt}$\begin{matrix}
 \circ &\circ &\circ\\
 \circ &\circ &\circ \\
 \circ  &\  &\  \\

  \end{matrix}$
 \end{tabular}
&
&\hspace{-5pt}\begin{tabular}{llllll}
2 &2 &2 &1 &1 &1\\
 &2 &2 &1 &1 &1\\
&&2 &1 & &
 \end{tabular}
& & \\
\begin{tabular}{ll}
\fbox{$\begin{matrix}
  \circ &\circ &\circ\\
 \circ &\circ &\circ \\
\end{matrix}$}
&\hspace{-5pt}$\begin{matrix}
  \circ\\
\
 \end{matrix}$
\end{tabular}
&
&\hspace{-5pt}\begin{tabular}{lll}
2 &2 &1\\
 &2 &\\

 \end{tabular}
&\
&\ \\
\\
\end{tabular}

\noindent{\it   Second Proof of  Theorem \ref{thmodd}.}
 Examining Corteel and Savage's second bijection $\gamma$ from $A$ to $D\times E$,
we see that it maps an anti-lecture hall composition of $n$ in $A$ with the first part not exceeding
$2k-1$  to a pair $(\alpha, \beta)$ in $D\times E$ such that $\beta$ is an
anti-lecture hall composition in $E$ with the first part $\beta_1$ not exceeding $2k-2$ and the sum of parts of $\alpha$ and $\beta$ equals $n$. In other words,
 $\gamma$ is a bijection between $F_{2k-1}$ and $D\times Q_{2k-2}$. Together with Theorem \ref{odd3}, we are led to a bijection between $F_{2k-1}$ and $D \times R_k$.

On the other hand, there is   a combinatorial interpretation of the left hand side of (\ref{generr}) in terms of the
 Durfee dissection of a partition, given by Andrews \cite{and79}.
 We observe that technique of Andrews easily extends to Durfee rectangle dissection of a
 partition. In this way,
 we find that the generating function of partitions in $R_k(n)$ is given by
 \begin{equation}\label{DE}
\sum_{n=0}^{\infty}|R_k(n)|q^n
=\sum_{N_1\geq N_2 \geq\ldots\geq N_{k-1}\geq 0}\frac{q^{N_1^2+\ldots+ N_{k-1}^2+N_1+\ldots +N_{k-1}}}
{(q)_{N_1-N_2}\ldots(q)_{N_{k-2}-N_{k-1}}
(q)_{N_{k-1}}}.
\end{equation}

 Setting $a=1$ in the generalization of the Rogers-Ramanujan identity (\ref{generr}) gives
 \begin{equation*}
\sum_{N_1\geq N_2 \geq \cdots \geq N_{k-1}\geq 0}
\frac{q^{N_1^2+\cdots+N_{k-1}^2+
N_1+\cdots+N_{k-1}}}{(q)_{N_1-N_2}\ldots(q)_{N_{k-2}-N_{k-1}}
(q)_{N_{k-1}}}
=\frac{(q,q^{2k},q^{2k+1};q^{2k+1})_{\infty}}
{(q;q)_{\infty}}.
\end{equation*}
Hence the generating function of partitions in  $R_k(n)$ can be expressed as follows
\begin{equation}\label{DE2}
\sum_{n=0}^{\infty}|R_k(n)|q^n
=\frac{(q,q^{2k},q^{2k+1};q^{2k+1})_{\infty}}
{(q;q)_{\infty}}.
\end{equation}
By the bijection between  $F_{2k-1}(n)$
 and $D\times R_k(n)$ we conclude that
\begin{equation}\label{F}
\sum_{n=0}^{\infty}|
F_{2k-1}(n)|q^n=\frac{(-q;q)_{\infty}(q,q^{2k},q^{2k+1};q^{2k+1})_{\infty}}{(q;q)_{\infty}}.
\end{equation}
It is easy to see that the right hand side of the above identity is the generating function of overpartitions in
 $H_{2k+1}(n)$.
This completes the proof.
\qed

 \noindent{\bf Acknowledgments.}  This work was supported by  the 973
Project, the PCSIRT Project of the Ministry of Education,  and the National Science
Foundation of China.

\end{document}